\documentclass[11pt,amsfonts]{article}
\usepackage[utf8]{inputenc}
\usepackage{graphicx}
\usepackage{latexsym}
\usepackage{amssymb}
\usepackage{amsmath}
\usepackage{enumerate}
\usepackage{color}
\usepackage{layout}
\usepackage{eufrak}
\usepackage{float}
\usepackage{dsfont}

\newtheorem{prop}{Proposition}
\newtheorem{lemma}{Lemma}

\newtheorem{corollary}{Corollary}
\newtheorem{theorem}{Theorem}
\newtheorem{remark}{Remark}

\def\real{{\mathord{{\rm I\kern-2.8pt R}}}}        % Fake blackboard bold R.
\def\inte{{\mathord{{\rm I\kern-2.8pt N}}}}

\def\sZZ{{\rm Z\kern-2.8ptem{}Z}}

\def\z{{\mathchoice
  {\sZZ}
  {\sZZ}
  {\rm Z\kern-0.30em{}Z}
  {\rm Z\kern-0.25em{}Z} }}
\def\sQQ{{\kern 0.27em \vrule height1.45ex width0.03em depth0em
          \kern-0.30em \rm Q}}
\def\qu{{\mathchoice
    {\sQQ}
    {\sQQ}
  {\kern 0.225em \vrule height1.05ex width0.025em depth0em \kern-0.25em \rm Q}
  {\kern 0.180em \vrule height0.78ex width0.020em depth0em \kern-0.20em \rm Q}
        }}
\def\sCC{{\kern 0.27em \vrule height1.45ex width0.03em depth0em
          \kern-0.30em \rm C}}
\def\complex{{\mathchoice
    {\sCC}
    {\sCC}
  {\kern 0.225em \vrule height1.05ex width0.025em depth0em \kern-0.25em \rm C}
  {\kern 0.180em \vrule height0.78ex width0.020em depth0em \kern-0.20em \rm C}
        }}

\newcommand{\E}{\mathbb{E}}

\newcommand{\ba}{\begin{array}}
\newcommand{\ea}{\end{array}}
\newcommand{\be}{\begin{equation}}
\newcommand{\ee}{\end{equation}}
\newcommand{\bea}{\begin{eqnarray}}
\newcommand{\eea}{\end{eqnarray}}
\newcommand{\beaa}{\begin{eqnarray*}}
\newcommand{\eeaa}{\end{eqnarray*}}

\def\z{\zeta}

\font\tenmath=msbm10 \font\sevenmath=msbm7 \font\fivemath=msbm5
\newfam\mathfam \textfont\mathfam=\tenmath
\scriptfont\mathfam=\sevenmath \scriptscriptfont\mathfam=\fivemath

\def \={{\buildrel {\rm (law)} \over =}}

\def\qed{ \hfill \vrule width.25cm height.25cm depth0cm\smallskip}

\newcommand{\basa}{\begin{assumption}}
\newcommand{\easa}{\end{assumption}}

\newcommand{\bas}{\begin{assum}}
\newcommand{\eas}{\end{assum}}

\newcommand{\ignore}[1]{}
\begin{document}

\renewcommand{\thefootnote}{\fnsymbol{footnote}}

\renewcommand{\thefootnote}{\fnsymbol{footnote}}

\title{On quadratic variations of the fractional-white wave equation}
\author{Radomyra Shevchenko $^{\ast}$ \vspace*{0.2in} \\
 $^{\ast}$ Fakult\"at f\"ur Mathematik,
LSIV, TU Dortmund \\
 Vogelpothsweg 87, 44227 Dortmund, Germany. \\
\quad radomyra.shevchenko@tu-dortmund.de \vspace*{0.1in}}

\maketitle

\begin{abstract} 
This paper studies the behaviour of quadratic variations of a stochastic wave equation driven by a noise that is white in space and fractional in time. Complementing the analysis of quadratic variations in the space component carried out in \cite{KT} and \cite{SST}, it focuses on the time component of the solution process. For different values of the Hurst parameter a central and a noncentral limit theorems are proved, allowing to construct parameter estimators and compare them to the finding in the space-dependent case. Finally, rectangular quadratic variations in the white noise case are studied and a central limit theorem is demonstrated.
\end{abstract}

\vskip0.3cm

{\bf 2010 AMS Classification Numbers:} Primary 60G15, 60G22, 60H15, 62F12; Secondary 62M30, 60F05.

\vskip0.3cm


\section{Introduction}

Statistical inference for stochastic partial differential equations (SPDEs) is an important and rapidly advancing branch of mathematical statistics. Usually under the framework of a Brownian field driving the equations new areas of applications are emerging (see e.g. \cite{PFABS} or \cite{ABJR}) and new methods are being developed for estimating the drift and volatility parameters in various settings (see \cite{C} for an extensive survey on the development of the subject in recent years).

One of the classical ideas for parameter estimation consists in considering so-called empirical power variations, that is, sums over increments of the solution process (either in the time or in the space component) raised to some power, see for instance \cite{BiT} or \cite{CD}. In particular, in a recent work \cite{HT} an in-depth study of quadratic variations for solutions of parabolic SPDEs is conducted on a space-time grid. 

In this context the development of stochastic calculus with respect to the fractional Brownian motion has led naturally to statistical inference for SPDEs driven by fractional noise either in the time or in the space component. Many authors have already investigated this topic over the last few decades (see, for example, \cite{BT} and \cite{DT}). For such equations the method of power variations can also be used in order to estimate the corresponding Hurst parameter of the driving noise analogously to the classical results for fractional Brownian motion and many associated processes (see the monograph \cite{T} for numerous examples).

In this paper we consider the stochastic wave equation with zero boundary conditions driven by a noise that is fractional in time and white in space. From the point of view of applications a solution to such an equation describes the motion of a randomly perturbed string. This equation and its properties has been described, for instance, in \cite{T} and \cite{BT}. In the paper \cite{KT} the authors study the behaviour of quadratic variations in the space coordinate for the Hurst parameter $H$ varying from $\frac{1}{2}$ to $\frac{3}{4}$ and in \cite{SST} the case $H>\frac{3}{4}$ is being considered. In both works the authors derive and analyse estimators for $H$.

The papers \cite{KT} and \cite{SST} have served as the starting point for the present manuscript. We study the behaviour of quadratic variations in the \textit {time} component of the wave equation solution. More precisely, if $u(t,x), t\geq 0, x\in \mathbb{R}$, denotes the solution to the wave equation with fractional-white noise, we consider the sequence of the centred (empirical) quadratic variations defined by
\begin{equation}\label{vn-intro}
V_N:=\frac{1}{N}\sum_{i=0}^{N-1}\left(u\left(\frac{i+1}{N},\,x\right)-u\left(\frac{i}{N},\,x\right) \right)^2-\mathbb E\left[\left(u\left(\frac{i+1}{N},\,x\right)-u\left(\frac{i}{N},\,x\right) \right)^2\right].
\end{equation}
We retrieve a standard threshold for processes in the fractional Brownian context and prove for the sequence $V_N$ a (quantitative) central limit theorem for the Hurst parameter $H$ between $\frac{1}{2}$ and $\frac{3}{4}$ as well as a noncentral limit theorem for $H$ above $\frac{3}{4}$, although the limiting object is different from the one obtained in \cite{SST} for space-dependent quadratic variations. Using these results and assuming that the mild solution $u$ is observed at discrete times and  at a fixed space location $x$,  we construct an estimator of  the parameter H from the observations $u(\frac{i}{N}, x)$ for $i=1,\dots ,N$. Based on the behaviour of the sequence  (\ref{vn-intro}), we prove that the estimator for $H$ is strongly consistent and asymptotically normal. Subsequently, we briefly compare this estimator to its space-dependent analogue from \cite{KT}. Furthermore, we introduce drift and volatility parameters into the equation and propose strongly consistent and asymptotically normal estimators for those. Finally, in the simpler scenario of Brownian noise (that is, for $H=\frac{1}{2}$) we consider rectangular, i.e. joint space-time, quadratic variations and prove a quantitative central limit theorem in this case. This allows us to construct a drift parameter estimator based on space-time observations and assess its asymptotic properties.

Methodically the results in this paper boil down to a meticulous analysis of the covariance structure of the solution to the wave equation (which is of independent interest) as well as to the application of classical techniques from the Malliavin-Stein toolkit such as the celebrated fourth moment theorem or the study of the cumulants in order to demonstrate convergence in distribution.

The paper is structured as follows. In Section \ref{S2} we briefly describe the setting and in Section \ref{S3} we study the covariance structure of the solution process in time. In Section \ref{S4} the main theorems are proved, namely a central limit theorem for $H<\frac{3}{4}$ and a noncentral limit theorem for $H>\frac{3}{4}$. Sections \ref{S5} and \ref{EstimC} deal with estimation questions for different settings related to the wave equation. Finally, in Section \ref{S7} several results are collected concerning rectangular quadratic variations in the simple case $H=\frac{1}{2}$. The paper ends with a concise appendix containing basic results and definitions from Malliavin calculus.

%The development of the stochastic calculus with respect to the fractional Brownian motion yields naturally to the statistical inference for SPDEs driven by fractional noise. 
%Many authors have already investigated this topic over the last few decades. While many works concern the estimation of the drift parameter for stochastic equation driven by fractional Brownian motion (we refer, among many others, to \cite{AV}, \cite{HN}, \cite{KB}, \cite{TV}), statistical inference for the Hurst parameter in the SPDEs driven by fractional noise has been considered less and attracted an increased research recently. 
%In order to obtain the asymptotic behavior in distribution of the temporal quadratic variation, we compute explicitly the correlation structure of the solution in time.  
%In order to estimate  the joint increments of the process u, we need to compute explicitly the covariance structure of the solution in time. 
%Assuming that the mild solution $u$ is observed at discrete times and  at a fixed space location x,  we construct estimator of  the parameter H from the observations $u(\frac{i}{N}, x)$. Based on the behavior of the sequence  (\ref{vn-intro}), we prove that the estimator for $H$ is consistent and asymptotically normal. 

\section{Preliminaries}\label{S2}

In this chapter we introduce the fractional-white wave equation and its solution and present the basic definitions used in our work.

The object of our study will be the solution to the following stochastic wave equation
\begin{equation}
\left\{
\begin{array}{rcl}\label{systeme wave}
\frac{\partial^2 u}{\partial t^2}(t,x)&=&\frac{\partial^2 u}{\partial x^2}
(t,x)+\dot W^H(t,x),\;t\geq 0,\;x \in \mathbb{R},\nonumber\\
\noalign{\vskip 2mm}
u(0, x)&=& 0, \quad x \in \mathbb{R},\nonumber\\
\noalign{\vskip 2mm} \frac{\partial u}{\partial t}(0,x) &=& 0,\quad
x \in \mathbb{R},
\end{array} \right.
\end{equation}
where $W^{H}$ is a fractional-white Gaussian noise which is defined as a real valued centred Gaussian field $W^{H}=\{ W^{H} (t,x); t \in  [0, T ], x\in \mathbb{R}\}$, over a given complete filtered probability space
$(\Omega, \mathcal{F}, ( \mathcal{F}_{t} )_{t\geqslant0} , \mathbb{P})$, with covariance function given by 
\begin{equation}\label{covW}
\mathbb{E}\left(W^{H} (t,x)W^ {H} (s,y)\right) = R_{H} (t, s)\min (x,y), \forall x,y \in \mathbb{R},
\end{equation}
where $R_{H}$ is the covariance of the fractional Brownian motion
 \begin{equation*}
R_{H}(t,s)=\frac{1}{2}\left(t^{2H}+ s^{2H}- \vert t-s \vert^{2H} \right), \hskip0.3cm s,t\geq 0.
\end{equation*}
We will assume throughout this work $H\in \left[\frac{1}{2}, 1\right).$ 

The solution of the equation (\ref{systeme wave})  is understood in the mild sense, that is, it is defined as a square-integrable
centered field $u=\left( u(t,x);\;t\in[0,T], x \in\mathbb{R} \right)$
 defined by
\begin{equation}\label{sol-wave-WFrac}
u(t,x)=\int_{0}^{t}\int_{\mathbb{R}}G_1(t-s,x-y)W^H(\mathrm{d}s,\mathrm{d}y),\quad t\geq 0, x\in \mathbb{R},
\end{equation}
where $G_{1}$ is the fundamental solution to the wave equation and the integral in (\ref{sol-wave-WFrac}) is a Wiener integral with respect to the Gaussian process $W^{H}$, that is, we have simply
\begin{equation}
\label{g1}
G_1(t,x)=\frac{1}{2}\mathds{1}_{\{|x|<t\}}.
\end{equation}

In the course of the paper we use the symbol $\sim$ to denote asymptotic equality (i.e. the ratio is tending to one), the symbol $\sim_c$ to denote asymptotic equality up to a constant, and the symbol $\lesssim$ to denote that the left side is asymptotically less or equal to the right side up to a constant (i.e. the ratio is asymptotically bounded by a constant).

\section{The temporal covariance structure}\label{S3}
The main factor in understanding the behaviour of a Gaussian process is determining its covariance structure which is calculated in this section.
\begin{theorem}
	For $\frac{1}{2}\leq H<\frac{3}{4}$he solution process $(u(t,\,x))_{t\in\mathbb R^+}$ for a fixed $x\in\mathbb R$ has the covariance structure
	\[\frac{4}{H}\mathbb E[u(t,\,x)u(s,\,x)]= \frac{1}{H(2H+1)}(t^{2H+1}+s^{2H+1})-\frac{2}{2H}t(t-s)^{2H}+\frac{2}{2H+1}(t-s)^{2H+1} \text{ for }t\geq s.\]
\end{theorem}
\begin{proof}
The proof for $H=\frac{1}{2}$ is given in Lemma \ref{CovRect} in Section 7 of this article.
For $H>\frac{1}{2}$ recall first that
	\[u(t,\,x)=\int_0^t\int_\mathbb R G_1(t-r,\,x-y)W^H(dr,\,dy)=\frac{1}{2}\int_0^t \int_\mathbb R 1_{\{|x-y|<t-r\}}W^H(dr,\, dy).\]
Using the isometry property we have with $\alpha_H = H(2H-1)$
\begin{eqnarray*}
\mathbb E [u(t,\,x)u(s,\,x)]&=&\alpha_H \int_0^t\int_0^s\int_\mathbb R  G_1(t-u,\,x-y) G_1(s-v,\,x-y)|u-v|^{2H-2}dydvdu\\
&=&\frac{\alpha_H}{4}\int_0^t\int_0^s\int_\mathbb R 1_{\{|x-y|<t-u\}}1_{\{|x-y|<s-v\}}|u-v|^{2H-2}dydvdu\\
&=&\frac{\alpha_H}{4}\int_0^t\int_0^s|u-v|^{2H-2}\int_{\mathbb R}1_{\{|x-y|<\min (t-u,\,s-v)\}}dydvdu.
\end{eqnarray*}
By direct computation we obtain
\[\int_{\mathbb R}1_{\{|x-y|<\min (t-u,\,s-v)\}}dy = 2\min (t-u,\,s-v).\]
Consequently,
\begin{eqnarray*}
	&\frac{4}{\alpha_H}&\mathbb E [u(t,\,x)u(s,\,x)]=\int_0^t\int_0^s2\min (t-u,\,s-v)|u-v|^{2H-2}dvdu\\
	&=&\int_0^t\int_0^s2\min (u,\,v)|t-u-s+v|^{2H-2}dvdu\\
	&=&\int_0^s\int_0^s2\min (u,\,v)|t-u-s+v|^{2H-2}dvdu+\int_s^t\int_0^s2\min (u,\,v)|t-u-s+v|^{2H-2}dvdu\\
	&=&\int_0^s\int_0^u2v|t-u-s+v|^{2H-2}dvdu+\int_0^s\int_u^s2u|t-u-s+v|^{2H-2}dvdu\\
	&&\qquad+\int_s^t\int_0^s2v|t-u-s+v|^{2H-2}dvdu\\
	&=&I_1+I_2+I_3.
\end{eqnarray*}
Let us assume $2s-t<0$ and analyse the three summands separately.
\begin{eqnarray*}
	I_1&=&\int_0^s 2v \int_v^s |t-u-s+v|^{2H-2}dudv\\
	&=&\int_0^s 2v \int_{s-t+v}^{2s-t} |u-v|^{2H-2}dudv\\
	&=&-\frac{1}{2H-1}\int_0^s2v((v-2s+t)^{2H-1}-(t-s)^{2H-1})dv\\
	&=&-\frac{1}{2H-1}\int_0^s 2v(v-2s+t)^{2H-1}dv+\frac{1}{2H-1}(t-s)^{2H-1}\int_0^s 2v dv\\
	&=&-\frac{1}{2H-1}\int_{-2s+t}^{t-s}2(v+2s-t)v^{2H-1}dv+\frac{1}{2H-1}(t-s)^{2H-1}s^2\\
	&=&\frac{1}{2H-1}\left((t-s)^{2H-1}s^2-\int_{-2s+t}^{t-s}2(v+2s-t)v^{2H-1}dv\right).
\end{eqnarray*}
For the second summand we obtain
\begin{eqnarray*}
	I_2&=& \int_0^s 2u\int_u^s |t-u-s+v|^{2H-2}dvdu\\
	&=&\int_0^s 2u\int_{t-s+u}^t |u-v|^{2H-2}dvdu\\
	&=&\frac{1}{2H-1} \int_0^s 2u ((t-u)^{2H-1}-(t-s)^{2H-1})du\\
	&=&\frac{1}{2H-1}\int_{t-s}^t 2(t-u)u^{2H-1}du-\frac{1}{2H-1}(t-s)^{2H-1}s^2\\
	&=&\frac{1}{2H-1}\left(\frac{2}{2H}t(t^{2H}-(t-s)^{2H})-\frac{2}{2H+1}(t^{2H+1}-(t-s)^{2H+1})-(t-s)^{2H-1}s^2\right).
\end{eqnarray*}
Finally, for the third summand we have
\begin{eqnarray*}
	I_3&=& \int_0^s 2v \int_s^t |t-u-s+v|^{2H-2}dudv\\
	&=&\int_0^s 2v \int_{2s-t}^s |u-v|^{2H-2}dudv\\
	&=&\int_0^s 2v \left (\int_{2s-t}^v (v-u)^{2H-2}du+\int_{v}^s (u-v)^{2H-2}du\right)dv\\
	&=&\frac{1}{2H-1}\int_0^s 2v((v-2s+t)^{2H-1}+(s-v)^{2H-1})dv\\
	&=&\frac{1}{2H-1}\int_{-2s+t}^{t-s}2(v+2s-t)v^{2H-1}dv + \frac{1}{2H-1}\int_0^s 2(s-v)v^{2H-1}dv\\
	&=&\frac{1}{2H-1}\left(\int_{-2s+t}^{t-s}2(v+2s-t)v^{2H-1}dv+\frac{2}{2H}s^{2H+1}-\frac{2}{2H+1}s^{2H+1}\right).
\end{eqnarray*}
Adding up the summands we obtain the result.\\
Let us turn to the case $2s-t\geq 0$. The first summand is
\begin{eqnarray*}
	I_1&=&\int_0^s 2v \int_{s-t+v}^{2s-t} |u-v|^{2H-2}dudv\\
	&=&\int_0^{2s-t} 2v \int_{s-t+v}^{2s-t} (u-v)^{2H-2}dudv+\int_{2s-t}^s 2v \int_{s-t+v}^{2s-t} (v-u)^{2H-2}dudv\\
	&=&\frac{1}{2H-1}\bigg(\int_0^{2s-t}2v(t-s)^{2H-1}dv+\int_0^{2s-t}2v(2s-t-v)^{2H-1}dv\\
	&&\qquad\qquad -\int_{2s-t}^s 2v(v-2s+t)^{2H-1}dv+\int_{2s-t}^s2v(t-s)^{2H-1}dv \bigg)\\
	&=&\frac{1}{2H-1} \left(\int_0^s 2vdv (t-s)^{2H-1}+\int_0^{2s-t}2v(2s-t-v)^{2H-1}dv-\int_{2s-t}^s 2v(v-2s+t)^{2H-1}dv\right)\\
	&=&\frac{1}{2H-1}\left((t-s)^{2H-1}s^2+\int_0^{2s-t}2v(2s-t-v)^{2H-1}dv-\int_{2s-t}^s 2v(v-2s+t)^{2H-1}dv\right).
\end{eqnarray*}
The second summand is the same as above, and for the third summand we obtain
\begin{eqnarray*}
	I_3&=&\int_0^s 2v \int_{2s-t}^s |u-v|^{2H-2}dudv\\
	&=&\int_0^{2s-t} 2v \int_{2s-t}^s (u-v)^{2H-2}dudv+ \int_{2s-t}^s 2v \int_{2s-t}^v (v-u)^{2H-2}dudv+\int_{2s-t}^s 2v \int_{v}^s (u-v)^{2H-2}dudv\\
	&=&\frac{1}{2H-1}\bigg(\int_0^{2s-t}2v(s-v)^{2H-1}dv-\int_0^{2s-t}2v(2s-t-v)^{2H-1}dv\\
	&&\qquad \qquad + \int_{2s-t}^s 2v(v-2s+t)^{2H-1}dv+\int_{2s-t}^s2v(s-v)^{2H-1}dv\bigg)\\
	&=&\frac{1}{2H-1}\left(\int_0^s 2(s-v)v^{2H-1}dv-\int_0^{2s-t}2v(2s-t-v)^{2H-1}dv+ \int_{2s-t}^s 2v(v-2s+t)^{2H-1}dv\right)\\
	&=&\frac{1}{2H-1}\left(\frac{2}{2H}s^{2H+1}-\frac{2}{2H+1}s^{2H+1}-\int_0^{2s-t}2v(2s-t-v)^{2H-1}dv+ \int_{2s-t}^s 2v(v-2s+t)^{2H-1}dv\right).
\end{eqnarray*}
Summing up $I_1$, $I_2$ and $I_3$ yields the same result.
\end{proof}

There are several remarks to be made concerning this result. First, the covariance is independent of space. Moreover, since the solution is Gaussian it follows directly from the covariance formula that it is a self-similar process in time. It can also be concluded from the formula that the process $u(\cdot,x)$ has a version with continuous paths with H\"older index below $H$, since
\[\mathbb E [(u(t,x)-u(s,x))^2]\lesssim |t-s|^{2H} \text{ for }t,s\geq 0,\]
and by Gaussianity
\[\mathbb E [(u(t,x)-u(s,x))^{2m}]\lesssim |t-s|^{2Hm} \text{ for }t,s\geq 0,\]
for $m\in\mathbb N$. The statement now follows by Kolmogorov's continuity criterion.\\

Next statements are concerned with the asymptotics of the covariance.

\begin{remark}
	In particular, we obtain for covariance of the increments:
	\begin{eqnarray*}
		&&\frac{4}{H}\mathbb E \left[\left(u\left(\frac{i}{N},\,x\right)-u\left(\frac{i-1}{N},\,x\right)\right)\left(u\left(\frac{j}{N},\,x\right)-u\left(\frac{j-1}{N},\,x\right)\right)\right]\\
		&&=\frac{1}{N^{2H+1}}\Big(\frac{2}{2H}\left(i(i-j+1)^{2H}-i(i-j)^{2H}+(i-1)(i-j-1)^{2H}-(i-1)(i-j)^{2H}\right)\\
		&&\qquad -\frac{2}{2H+1}\left((i-j-1)^{2H+1}-2(i-j)^{2H+1}+(i-j+1)^{2H+1}\right)\Big)
	\end{eqnarray*}
if $i>j$ and
		\begin{eqnarray*}
		&&\frac{4}{H}\mathbb E \left[\left(u\left(\frac{i+1}{N},\,x\right)-u\left(\frac{i}{N},\,x\right)\right)^2\right]=2\frac{1}{N^{2H+1}}\left(\frac{i}{H}+\frac{1}{H(2H+1)}\right).
	\end{eqnarray*}
\end{remark}

%\begin{lemma}
%	For $i>j$ we have asymptotically (for large $(i-j)$)
%	\begin{eqnarray*}
%		&&i(i-j+1)^{2H}-i(i-j)^{2H}+(i-1)(i-j-1)^{2H}-(i-1)(i-j)^{2H}\\
%		&&\quad\sim 2H(i-j)^{2H-1}+2H(2H-1)(i-j)^{2H-2}(i-1)
%	\end{eqnarray*}
%as well as
%	\[(i-j-1)^{2H+1}-2(i-j)^{2H+1}+(i-j+1)^{2H+1}\sim2H(2H+1)(i-j)^{2H-1}.\]
%\end{lemma}
%\begin{proof}
%	For the first statement we calculate with binomial expansion
%	\begin{eqnarray*}
%	&&i(i-j+1)^{2H}-i(i-j)^{2H}+(i-1)(i-j-1)^{2H}-(i-1)(i-j)^{2H}\\
%	&&\quad = (i-1)(i-j)^{2H}\left(\frac{i}{i-1}\left(1+\frac{1}{i-j}\right)^{2H}-\frac{i}{i-1}+\left(1-\frac{1}{i-j}\right)^{2H}-1\right)\\
%	&&\quad =(i-1)(i-j)^{2H}\Big(1+\frac{2H}{i-j}+\frac{H(2H-1)}{(i-j)^2}+O\left(\frac{1}{(i-j)^3}\right)-2\\
%	&&\qquad+\frac{1}{i-1}\left(1+\frac{2H}{i-j}+\frac{H(2H-1)}{(i-j)^2} + O\left(\frac{1}{(i-j)^3}\right)\right)-\frac{1}{i-1}\\
%	&&\qquad +1-\frac{2H}{i-j}+\frac{H(2H-1)}{(i-j)^2}+O\left(\frac{1}{(i-j)^3}\right)\Big)\\
%	&&\quad\sim 2H(i-j)^{2H-1}+2H(2H-1)(i-j)^{2H-2}(i-1).
%	\end{eqnarray*}
%The second statement is well known and can be found, for instance, in \cite{coeur}.
%\end{proof}
\begin{corollary}\label{CorCov}
	Note that for $i>j$ we can write the covariance function as follows:
\begin{eqnarray*}
&& \mathbb E \left[\left(u\left(\frac{i}{N},\,x\right)-u\left(\frac{i-1}{N},\,x\right)\right)\left(u\left(\frac{j}{N},\,x\right)-u\left(\frac{j-1}{N},\,x\right)\right)\right]\\
&&\qquad\qquad=\frac{H}{4N^{2H+1}}\left(\psi_1(i-j)+i\psi_2(i-j)\right),
\end{eqnarray*}
	where
	\[\psi_1(k)=\frac{2}{2H}\left(k^{2H}-(k-1)^{2H}\right)-\frac{2}{2H+1}\left((k+1)^{2H+1}-2k^{2H+1}+(k-1)^{2H+1}\right)\]
	and
	\[\psi_2(k)=\frac{2}{2H}\left((k+1)^{2H}-2k^{2H}+(k-1)^{2H}\right)\]
	with the following asymptotics for large $k$:
	\[\psi_1(k)\sim 2(1-2H)k^{2H-1},\qquad \psi_2(k)\sim 2(2H-1)k^{2H-2}.\]
These expressions are obtained using the binomial expansion applied for $\left(1 + \frac{1}{k}\right)^{2H}$,  $\left(1 - \frac{1}{k}\right)^{2H}$, $\left(1 + \frac{1}{k}\right)^{2H+1}$ and $\left(1 - \frac{1}{k}\right)^{2H+1}$.\\
	We, moreover, obtain for large $i-j$ using the same asymptotics:
	\[\psi_1(i-j)+i\psi_2(i-j)\sim (2H-1)(2j-1)(i-j)^{2H-2}.\]
\end{corollary}

\section{The temporal quadratic variations}\label{S4}
For the solution $u$ of the wave equation we define its quadratic variation in time,
\[V_N:=\frac{1}{N}\sum_{i=0}^{N-1}\left(u\left(\frac{i+1}{N},\,x\right)-u\left(\frac{i}{N},\,x\right) \right)^2-\mathbb E\left[\left(u\left(\frac{i+1}{N},\,x\right)-u\left(\frac{i}{N},\,x\right) \right)^2\right].\]
For simplicity let us denote $u_i:=u\left(\frac{i}{N},\,x\right)$ for some fixed $x\in\mathbb R$.

\subsection{Renormalization of $V_N$}
\begin{prop}\label{EVn}
As $N$ tends to infinity, we have asymptotically $\mathbb E V_N^2\sim_c N^{-4H-1}$ for $\frac{1}{2}\leq H< \frac{3}{4}$ and $\mathbb E V_N^2\sim_c N^{-4}$ for $\frac{3}{4}<H<1$ up to some constants made exact in the proof.
\end{prop}
\begin{proof}
	We have by reordering the sum and putting together the non-diagonal summands that appear twice
	\begin{eqnarray*}
		\mathbb E[V_N^2]&=&\frac{1}{N^2}\sum_{i=0}^{N-1}\sum_{j=0}^{N-1}2\mathbb E [(u_{i+1}-u_i)(u_{j+1}-u_j)]^2\\
		&=&\frac{4}{N^2}\sum_{j=0}^{N-1}\sum_{i=j+1}^{N-1}\mathbb E [(u_{i+1}-u_i)(u_{j+1}-u_j)]^2+\frac{2}{N^2}\sum_{i=0}^{N-1}\mathbb E [(u_{i+1}-u_i)^2]^2=:v_{N,1}+v_{N,2}.
	\end{eqnarray*}
The non-diagonal summands with $|i-j|$ less than a certain constant are at most of order $N^{-4H-2}$ and can therefore be ignored in the asymptotics up to constants. We obtain
\begin{eqnarray*}
	v_{N,1}&\sim&\frac{4}{N^2}\frac{H^2}{16 N^{4H+2}}\sum_{j=0}^{N-1}\sum_{i=j+1}^{N-1} (2(2H-1)j(i-j)^{2H-2})^2\\
	&=&\frac{16(2H-1)^2}{N^2}\frac{H^2}{16 N^{4H+2}}\sum_{j=0}^{N-1}j^2\sum_{k=1}^{N-j-1}k^{4H-4}.
\end{eqnarray*} 
If $H>\frac{3}{4}$, $k^{4H-4}$ is not summable and $v_{N,1}$ is asymptotically equal to
\begin{eqnarray*}
	&&\frac{(2H-1)^2}{(4H-3)N^2}\frac{H^2}{N^{4H+2}}\sum_{j=0}^{N-1}j^2(N-j-1)^{4H-3}\\
	&&\quad =\frac{(2H-1)^2}{(4H-3)N^2} \frac{H^2}{N^{4H+2}}\sum_{j=0}^{N-1}(N-j-1)^2j^{4H-3}\\
	&&\quad\sim \frac{2(2H-1)^2}{(4H-3)\dots (4H) N^2} \frac{H^2}{N^{4H+2}} N^{4H}\\
	&&\quad = \frac{H(2H-1)}{4(4H-1)(4H-3)}N^{-4}.
\end{eqnarray*}

If $H<\frac{3}{4}$, $k^{4H-4}$ is summable. To obtain the precise constant we recall that
\begin{eqnarray*}
	v_{N,1}&=&\frac{4}{N^2}\frac{H^2}{16 N^{4H+2}}\sum_{i=0}^{N-1}\sum_{j=0}^{i-1} (\psi_1(i-j)+i\psi_2(i-j))^2\\
	&=&\frac{4}{N^2}\frac{H^2}{16 N^{4H+2}}\sum_{i=0}^{N-1}\sum_{j=0}^{i-1}(\psi_1(i-j)^2+2i\psi_2(i-j)\psi_1(i-j)+i^2\psi_2(i-j)^2).
\end{eqnarray*}
One can easily see with Corollary \ref{CorCov} that the first two summands are of order $N^{-4}$ while the third one is of order $N^{-4H-1}$ and dominates the other two. Therefore, we have
\begin{eqnarray*}
&&\lim_{N\to\infty} N^{4H+1}v_{N,1}=\lim_{N\to\infty} N^{4H+1}\frac{4}{N^2}\frac{H^2}{16 N^{4H+2}}\sum_{i=0}^{N-1}i^2 \sum_{j=0}^{i-1} \psi_2(i-j)^2\\
&&\qquad =\lim_{N\to\infty} \frac{H^2}{4}N^{-3}\sum_{j=1}^N \psi_2(j)^2\sum_{i=j}^N i^2=\frac{H^2}{12}\sum_{j=1}^\infty \psi_2(j)^2,
\end{eqnarray*}
which is summable.

Finally, for the diagonal we calculate
\begin{eqnarray*}
	v_{N,2}&=&\frac{2}{N^2}\sum_{i=0}^{N-1}\left(\frac{H}{2}\frac{1}{N^{2H+1}}\left(\frac{i}{H}+\frac{1}{H(2H+1)}\right)\right)^2\\
	&=&\frac{H^2}{2N^{4H+4}}\sum_{i=0}^{N-1}\left(\frac{i}{H}+\frac{1}{H(2H+1)}\right)^2\\
	&\sim &\frac{H^2}{2N^{4H+4}}\sum_{i=0}^{N-1}\left(\frac{i}{H}\right)^2\sim\frac{1}{6}N^{-4H-1}.\\
\end{eqnarray*}
For $H>\frac{3}{4}$ the term $N^{-4}$ is slower than $N^{-4H-1}$, and the claim follows with nonzero limiting constants.\\
More precisely, for $H<\frac{3}{4}$ we obtain
$$\lim_{N\to\infty} N^{4H+1}\mathbb E V_N^2=\frac{H^2}{12}\sum_{j=1}^\infty \psi_2(j)^2+\frac{1}{6}=\frac{H^2}{24}\sum_{j=-\infty}^\infty \psi_2(|j|)^2.$$
For $H>\frac{3}{4}$ we can write
\begin{eqnarray*}
	N^{4}\E(V_N^{2})&=&\frac{H^2}{4N^{4H}}\sum_{i=0}^{N-1}\sum_{j=0}^{i-1}\psi_{1}(i-j)^{2}+ \frac{H^2}{2N^{4H}} \sum_{i=0}^{N-1}\sum_{j=0}^{i-1} i \psi_{1}(i-j)\psi_{2}(i-j)\\
	&&+\frac{H^{2}}{4N^{4H}}\sum_{i=0}^{N-1}\sum_{j=0}^{i-1}i^{2}\psi_{2}(i-j) + N^{4} v_{N,2}\\
	\lim_{N \to \infty}N^{4}\E(V_N^2) &=& \lim_{N \to \infty}  \frac{H^2}{4N^{4H}}\sum_{j=0}^{N-1}\sum_{i=j+1}^{N-1}\psi_{1}(j)^{2}+ \frac{H^2}{2N^{4H}} \sum_{j=0}^{N-1}\psi_{1}(j)\psi_{2}(j)\sum_{i=j+1}^{N-1} i \\
	&& +\frac{H^{2}}{4N^{4H}}\sum_{j=0}^{N-1}\psi_{2}(j)^2\sum_{i=j+1}^{N-1}i^{2},\\
%	& = & \frac{4H^{2}(1-2H)^2}{4H-1}- \frac{H^2(1-2H)^2}{4H-2} + \frac{H^2(2H-1)^2}{3(4H-3)}\\
%	&=& \frac{H^{2}(1-H)^{2}(16H^2-24H+11)}{3(4H-1)(4H-3)(4H-2)}
\end{eqnarray*}
which defines the normalising constant.
\end{proof}

\begin{remark}\label{sigma2}
Note that with the notation from \cite{KT}, namely
$$\varphi_H(k)=\frac{1}{2}\left((k+1)^{2H}-2k^{2H}+(k-1)^{2H}\right),$$
the precise limiting constant for $H<\frac{3}{4}$ equals
$$\sigma^2=\frac{1}{6}\sum_{k\in\mathbb Z}\varphi_H(k)^2.$$
\end{remark}

Now we know which normalisation is needed to prove limit theorems. We consider $F_N:=\frac{V_N}{\sqrt{\mathbb E [V_N^2]}}$. For $H<\frac{3}{4}$ we have $F_N\sim_cN^{2H+\frac{1}{2}}V_N$ pointwise.

\subsection{Central limit theorem and rate of convergence}\label{CLT}
To establish the central limit theorem of the quadratic variations, we will use tools from the Malliavin-Stein framework. A short introduction of the necessary terminology and classical identities can be found in the Appendix. The principal statement necessary for the proof of the theorem is Theorem 5.2.6 in \cite{N-P}, which is a version of the fourth moment theorem. For convenience of the reader we recall it in the following.

\begin{theorem}\label{clt-fract-white}
Fix $q\geq 1$. Let $(G_{N})_{N\geq 1} = (I_{q}(g_{N}))_{N\geq 1}$
with $g_{N} \in \mathcal{H}^{\odot q}$, be a sequence of random
variables belonging to the $q$th Wiener chaos such that
$$\mathbb{E}(G_{N}^{2})\stackrel{N\to \infty}{\to} s ^{2} >0.$$
Then $G_{N} $ converges in law to $Z\sim\mathcal{N}(0,1)$ if and only if
$$\Vert DG_{N} \Vert ^{2}_{\mathcal{H}}\stackrel{N\to \infty}{\to}  qs^{2}.$$
Furthermore,
\begin{equation*}
d\left(G_{N};\mathcal{N}(0,1) \right) \leq C \sqrt{\mathbf{Var}
\Big(\frac{1}{q}{\Vert DG_{N} \Vert}^{2} _{\mathcal{H}} }\Big),
\end{equation*}
where $d$ is either the distance of Kolmogorov, the distance in 
Total Variation or the Wasserstein distance. 
\end{theorem}

From now on, fix $x \in \mathbb{R}$  and denote by $\mathcal{H}$ the Hilbert space associated to the Gaussian solution process $(u(t, x))_{t \in [0,1]}$. 
This Hilbert space is defined as the closure of the set of indicator functions $\{1_{[0,t]}, t \in [0,1]\}$ with respect to the inner product

$$\langle 1_{[0,t]}, 1_{[0,s]}\rangle_{\mathcal{H}} = \E\left(u(t, x)u(s, x)\right), \mbox{for a fixed} \; x \in \mathbb{R}.$$

Denoting by $I_q(\cdot)$ the $q$th multiple integral with respect to $(u(t, x))_{t \in [0,1]}$, we can write
$$V_N=I_2\left(\frac{1}{N}\sum_{i=0}^{N-1}1_{\left[\frac{i}{N}, \frac{i+1}{N}\right]}^{\otimes 2}\right)$$
using the product rule \eqref{prod}. Now we can formulate the central limit theorem for $F_N$, the normalised version of $V_N$.

\begin{theorem}\label{clt-time}
For $\frac{1}{2}\leq H<\frac{3}{4}$ the sequence $F_N$ converges in law to $N(0,1)$ as $N$ tends to infinity. Moreover,
\[d(F_N,\,N(0,1))\lesssim
\begin{cases}
N^{-\frac{1}{2}} & \text{if $H \in [\frac{1}{2}, \frac{5}{8})$ }, \\
N^{-\frac{1}{2}} \log (N)^{\frac{3}{2}} & \text{if $H =\frac{5}{8}$}, \\
N^{4H-3} & \text{if $H \in (\frac{5}{8},\frac{3}{4})$}
\end{cases}\]
with $d$ either Kolmogorov, Wasserstein or total variation distance.
\end{theorem}
\begin{proof}
	We have
	\[D_sF_N = c N^{2H+\frac{1}{2}}\frac{1}{N}\sum_{i=0}^{N-1}I_1(1_{[\frac{i}{N},\frac{i+1}{N}]}(\cdotp))1_{[\frac{i}{N},\frac{i+1}{N}]}(s)\]
	with some constant $c$.
	Consequently, we obtain
	\begin{eqnarray*}
		&&\operatorname{Var}(\frac{1}{2}\|DF_N\|^2_{\mathcal H})\\
		&&\quad\sim_c \frac{1}{N^4}N^{8H+2}\sum_{ijkl=0}^{N-1}\langle 1_{[\frac{i}{N},\frac{i+1}{N}]},1_{[\frac{j}{N},\frac{j+1}{N}]} \rangle_{\mathcal H}\langle 1_{[\frac{k}{N},\frac{k+1}{N}]},1_{[\frac{j}{N},\frac{j+1}{N}]}\rangle_{\mathcal H}\\
&&\qquad\qquad\times  \langle 1_{[\frac{k}{N},\frac{k+1}{N}]},1_{[\frac{l}{N},\frac{l+1}{N}]}\rangle_{\mathcal H} \langle 1_{[\frac{i}{N},\frac{i+1}{N}]},1_{[\frac{l}{N},\frac{l+1}{N}]} \rangle_{\mathcal H}.
	\end{eqnarray*}
For the off-diagonal terms the values of the scalar products are determined by Corollary \ref{CorCov}. By symmetry we have to differentiate between five cases, namely summands with $i\psi_2(i-j)$ appearing $0$, $1$, $2$, $3$ or $4$ times. Let us consider those.

For zero appearances of $\psi_2$ and the sum going over non-diagonal terms we compute
\begin{eqnarray*}
	&&\frac{1}{N^4}N^{8H+2}\frac{1}{N^{8H+4}}\sum_{ijkl=0}^{N-1} \psi_1(|i-j|)\psi_1(|j-k|)\psi_1(|k-l|)\psi_1(|i-l|)\\
	&&=_cN^{-6} N \sum_{k=-N}^N (\psi_1\ast\psi_1)^2(|k|)\lesssim N^{-5}\left(\sum_{k=-N}^N |\psi_1(|k|)|^{\frac{4}{3}}\right)^3\\
	&&\lesssim N^{-5} \left(\sum_{k=-N}^N k^{(2H-1)\frac{4}{3}}\right)^3\lesssim N^{8H-6}
\end{eqnarray*}
by Young's inequality and Corollary \ref{CorCov}. For one appearance of $\psi_2$ we get is a similar manner (additionally using the Cauchy-Schwarz inequality)
\begin{eqnarray*}
	&&\frac{1}{N^4}N^{8H+2}\frac{1}{N^{8H+4}}\sum_{ijkl=0}^{N-1} \psi_1(|i-j|)\psi_1(|j-k|)\psi_1(|k-l|)\psi_2(|i-l|)\max(i,l)\\
	&&\lesssim N^{-4}\sum_{k=-N}^N (\psi_1\ast\psi_1)(|k|)(\psi_1\ast\psi_2)(|k|) \lesssim N^{-4} \|\psi_1\|^3_{{\frac{4}{3}}}\|\psi_2\|_{{\frac{4}{3}}}\\
	&&\lesssim N^{-4} \left(\sum_{k=-N}^N k^{(2H-1)\frac{4}{3}}\right)^{\frac{9}{4}}\left(\sum_{k=-N}^N k^{(2H-2)\frac{4}{3}}\right)^{\frac{3}{4}}\\
	&&\lesssim
	\begin{cases}
		N^{6H-\frac{19}{4}} & \text{if $H \in [\frac{1}{2}, \frac{5}{8})$ }, \\
		N^{6H-\frac{19}{4}} \log (N)^{\frac{3}{4}} & \text{if $H =\frac{5}{8}$}, \\
		N^{8H-6} & \text{if $H \in (\frac{5}{8},\frac{3}{4})$}.
	\end{cases}
\end{eqnarray*}
Similarly, for two appearances of $\psi_2$ we obtain
\begin{eqnarray*}
	&&\frac{1}{N^4}N^{8H+2}\frac{1}{N^{8H+4}}\sum_{ijkl=0}^{N-1} \psi_1(|i-j|)\psi_1(|j-k|)\psi_2(|k-l|)\max(k,l)\psi_2(|i-l|)\max(i,l)\\
	&&\lesssim N^{-3} \|\psi_1\|^2_{{\frac{4}{3}}}\|\psi_2\|^2_{{\frac{4}{3}}} \lesssim
	\begin{cases}
		N^{4H-\frac{7}{2}} & \text{if $H \in [\frac{1}{2}, \frac{5}{8})$ }, \\
		N^{4H-\frac{7}{2}} \log (N)^{\frac{3}{2}} & \text{if $H =\frac{5}{8}$}, \\
		N^{8H-6} & \text{if $H \in (\frac{5}{8},\frac{3}{4})$}.
	\end{cases}
\end{eqnarray*}
For three appearances of $\psi_2$ we get
\begin{eqnarray*}
	&&\frac{1}{N^4}N^{8H+2}\frac{1}{N^{8H+4}}\sum_{ijkl=0}^{N-1} \psi_1(|i-j|)\psi_2(|j-k|)\max(j,k)\psi_2(|k-l|)\max(k,l)\psi_2(|i-l|)\max(i,l)\\
	&&\lesssim N^{-2} \|\psi_1\|_{{\frac{4}{3}}}\|\psi_2\|^3_{{\frac{4}{3}}} \lesssim
	\begin{cases}
		N^{2H-\frac{9}{4}} & \text{if $H \in [\frac{1}{2}, \frac{5}{8})$ }, \\
		N^{2H-\frac{9}{4}} \log (N)^{\frac{3}{4}} & \text{if $H =\frac{5}{8}$}, \\
		N^{8H-6} & \text{if $H \in (\frac{5}{8},\frac{3}{4})$}.
	\end{cases}
\end{eqnarray*}
Finally, for four appearances of $\psi_2$ we obtain
\begin{eqnarray*}
	&&\frac{1}{N^4}N^{8H+2}\frac{1}{N^{8H+4}}\sum_{ijkl=0}^{N-1} \psi_2(|i-j|)\max(i,j)\psi_2(|j-k|)\max(j,k)\psi_2(|k-l|)\max(k,l)\psi_2(|i-l|)\max(i,l)\\
	&&\lesssim N^{-1} \|\psi_2\|^4_{{\frac{4}{3}}} \lesssim
	\begin{cases}
		N^{-1} & \text{if $H \in [\frac{1}{2}, \frac{5}{8})$ }, \\
		N^{-1} \log (N)^{3} & \text{if $H =\frac{5}{8}$}, \\
		N^{8H-6} & \text{if $H \in (\frac{5}{8},\frac{3}{4})$}.
	\end{cases}
\end{eqnarray*}
This last case showcases the slowest speeds of convergence, thus determining the behaviour of the non-diagonal part. Using the same technique, Young's and Cauchy-Schwarz inequalities, we obtain that also the diagonal terms converge faster or as fast as this last term. In total, we obtain
\begin{eqnarray*}
\operatorname{Var}(\frac{1}{2}\|DF_N\|^2_{\mathcal H})\lesssim
\begin{cases}
	N^{-1} & \text{if $H \in [\frac{1}{2}, \frac{5}{8})$ }, \\
	N^{-1} \log (N)^{3} & \text{if $H =\frac{5}{8}$}, \\
	N^{8H-6} & \text{if $H \in (\frac{5}{8},\frac{3}{4})$},
\end{cases}
\end{eqnarray*}
and the claim follows.
\end{proof}

\begin{remark}
	There seems to be a certain stability to the behaviour of quadratic variations of processes associated to fractional Brownian motion for $H<\frac{3}{4}$. The bounds on the speed of convergence obtained here copy those obtained for the space component and are the same as the bounds for the fractional Brownian motion itself.
\end{remark}

\subsection{Noncentral limit theorem}

Let us focus on the case $H>\frac{3}{4}$. Recall that in this case we have $\mathbb E [V_N^2]N^4\to k$ for some positive constant $k$ as $N$ tends to infinity.

\begin{theorem}\label{nclt}
For $H>\frac{3}{4}$ the sequence $\frac{V_N}{\sqrt{\mathbb E [V_N^2]}}$ converges in distribution to a random variable with cumulants
\begin{eqnarray*}
&&k_m=2^{2m-1}(2H-1)^m k^{-m\slash 2}(m-1)!\\
&&\qquad \times \int_{[0,1]^m}|x_1-x_2|^{2H-2}\min (x_1,x_2)|x_2-x_3|^{2H-2}\min (x_2,x_3)\\
&&\qquad\times\dots |x_m-x_1|^{2H-2}\min (x_m,x_1)dx_1\dots dx_m
\end{eqnarray*}
for $m\geq 3$ as $N$ tends to infinity.
\end{theorem}

\begin{proof}
By construction the random variables $\frac{V_N}{\sqrt{\mathbb E [V_N^2]}}$ are elements of the second chaos in the space generated by $W^H$, and hence in order to demonstrate convergence in distribution it is enough to determine the limits of all the cumulants of this sequence.

By the cumulants formula (see e.g. \cite{T}) for $G=I_2(f)$
\begin{equation*}
k_{m}(G)=  2^{m-1}(m-1)! \int_{ \mathbb{R}  ^{m} }d { u}_{1}\ldots d{ u}_{m} f({ u}_{1},
{ u}_{2}) f({ u}_{2}, { u}_{3})\ldots f({ u}_{m-1}, { u}_{m}) f({ u}_{m},{  u}_{1}),
\end{equation*}
following the same steps as in \cite{SST}, we obtain
\begin{eqnarray*}
&&k_m\left(\frac{V_N}{\sqrt{\mathbb E [V_N^2]}}\right)=2^{m-1}(m-1)!\left(\frac{1}{N\sqrt{\mathbb E [V_N^2]}}\right)^m\\
&&\qquad\times  \sum_{j_1,\dots , j_m=1}^N\mathbb E [(u_{j_1+1}-u_{j_1})(u_{j_2+1}-u_{j_2})]\dots \mathbb E [(u_{j_m+1}-u_{j_m})(u_{j_1+1}-u_{j_1})].
\end{eqnarray*}
By Corollary \ref{CorCov} and using the formula
\[|i-j+1|^{2H}-2|i-j|^{2H}+|i-j-1|^{2H}=2H(2H-1)\int_i^{i+1}\int_j^{j+1}|u-v|^{2H-2}dudv\]
we can write
\begin{eqnarray*}
&& \mathbb E [(u_{j_1+1}-u_{j_1})(u_{j_2+1}-u_{j_2})]=N^{-2H-1}\Big(\frac{2}{2H}|j_1-j_2|^{2H}\left(1-\left(1-\frac{1}{|j_1-j_2|}\right)^{2H}\right)\\
&&\qquad -\frac{2}{2H+1}2H (2H+1)|j_1-j_2|^{2H-1}\int_0^1\int_0^1 \left|1+\frac{u-v}{j_1-j_2}\right|^{2H-1}dudv\\
&&\qquad +\max (j_1,j_2)\frac{2}{2H}2H(2H-1)|j_1-j_2|^{2H-2}\int_0^1\int_0^1 \left|1+\frac{u-v}{j_1-j_2}\right|^{2H-2}dudv\Big).
\end{eqnarray*}
Using the formula for binomial series together with the dominated convergence theorem we can see that, as $N$ goes to infinity, this expression becomes equivalent to
\begin{eqnarray*}
&& N^{-2H-1}\Big(\frac{2}{2H}2H|j_1-j_2|^{2H-1} -4H|j_1-j_2|^{2H-1} +\max (j_1,j_2)2(2H-1)|j_1-j_2|^{2H-2}\Big)\\
&&\qquad = N^{-2}\left|\frac{j_1-j_2}{N}\right|^{2H-2}\frac{1}{N}((2-4H)|j_1-j_2|+\max (j_1,j_2)2(2H-1) )\\
&&\qquad=N^{-2}2(2H-1)\left|\frac{j_1-j_2}{N}\right|^{2H-2}\min \left(\frac{j_1}{N},\frac{j_2}{N}\right).
\end{eqnarray*}
Therefore, as $N$ tends to infinity, the $m$th cumulant is equivalent to
\begin{eqnarray*}
&& 2^{m-1}(m-1)N^{2m}!k^{-m\slash 2}N^{-m}N^{-2m} 2^m (2H-1)^m\sum_{j_1,\dots , j_m=1}^N \left|\frac{j_1-j_2}{N}\right|^{2H-2}\min \left(\frac{j_1}{N},\frac{j_2}{N}\right)\\
&&\qquad \times \left|\frac{j_m-j_1}{N}\right|^{2H-2}\min \left(\frac{j_m}{N},\frac{j_1}{N}\right).
\end{eqnarray*}
As $N$ tends to infinity, this expression converges to
\begin{eqnarray*}
&&2^{2m-1}(2H-1)^m k^{-m\slash 2}(m-1)!\\
&&\qquad \times \int_{[0,1]^m}|x_1-x_2|^{2H-2}\min (x_1,x_2)\dots |x_m-x_1|^{2H-2}\min (x_m,x_1)dx_1\dots dx_m.
\end{eqnarray*}
By Lemma 3.3 in \cite{BT} this integral is finite.\qed
\end{proof}

\begin{remark}
Unlike for the case $H<\frac{3}{4}$, the limiting objects for $H>\frac{3}{4}$ reflect the differences in the covariance structure with more detail. While the cumulants of the limit are similar to those of the Rosenblatt process (limiting object for the quadratic variations of the fractional Brownian motion) as well as to those of the corresponding result for the space component (see \cite{SST}), they are different from those two cases.\\

The case $H=\frac{3}{4}$ has not been considered in this paper for brevity reasons: A central limit theorem analogous to Theorem \ref{clt-time} is expected in this case with an additional factor $\log N$ in the speed of convergence, similarly to the result obtained in \cite{SST}.
\end{remark}
\section{Estimation of the Hurst parameter}\label{S5}
In this chapter we define and study an estimator for the Hurst parameter $H$ based on equidistant temporal observations on the interval $[0,1]$ of the solution $u$ at a fixed, possibly unknown point $x\in\mathbb R$ in space.\\

Let us define
$$S_N:=\frac{1}{N}\sum_{i=0}^{N-1}\left( u\left(\frac{i+1}{N},x\right)-u\left(\frac{i}{N},x\right)\right)^2$$
and note that $V_N=S_N-\mathbb E S_N$. We have
$$\mathbb E S_N=\frac{1}{4}N^{-2H}+\frac{1}{4}\frac{1-2H}{1+2H}N^{-2H-1},$$
and therefore,
$$\log (4\mathbb E S_N)=-2H\log N+\log \left(1+\frac{1-2H}{1+2H}N^{-1}\right)\sim -2H\log N + \frac{1-2H}{1+2H}N^{-1}.$$
Let us now define the estimator for $H$ as follows:
$$\hat H:=-\frac{\log 4+\log (S_N)}{2\log N}.$$
Then we obtain
\begin{eqnarray*}
-2\hat H &=&\frac{\log (\mathbb E S_N+V_N)+\log 4}{\log N}=\frac{\log (\mathbb E S_N(1+\frac{V_N}{\mathbb E S_N}))+\log 4}{\log N}\sim \frac{\frac{V_N}{\mathbb E S_N}+\log (4\mathbb E S_N)}{\log N}\\
&\sim & \frac{\frac{V_N}{\mathbb E S_N}-2H\log N + \frac{1-2H}{1+2H}N^{-1}}{\log N}.
\end{eqnarray*}
Consequently,
$$2(H-\hat H)\sim \frac{1}{\log N}\left(\frac{V_N}{\mathbb E S_N}+ \frac{1-2H}{1+2H}N^{-1}\right)\sim \frac{1}{\log N}\left(V_N N^{2H}+ \frac{1-2H}{1+2H}N^{-1}\right).$$
By Proposition \ref{EVn} and hypercontractivity \eqref{hyper}
$$\frac{V_N}{\mathbb E S_N}\sim V_N N^{2H}$$
goes to zero almost surely, hence ensuring strong consistency of the estimator.
Moreover, it follows directly from the CLT that for $\frac{1}{2}\leq H<\frac{3}{4}$
$$(H-\hat H)N\log N\stackrel{d}{\to} N\left(\frac{1-2H}{1+2H},\,\frac{\sigma ^2}{4}\right)$$
with $\sigma^2 = \frac{1}{6}\sum_{k\in\mathbb Z}\varphi_H(k)^2$ as in Remark \ref{sigma2}, showing that the estimator is asymptotically normal. For $H>\frac{3}{4}$ the term $(H-\hat H)N^{2-2H}\log N$ converges to a non-normal second chaos random variable as a consequence of Theorem \ref{nclt}.\\
The following theorem summarises our findings.
\begin{theorem}
Suppose that the solution field of the equation \eqref{systeme wave} is observed at some (possibly unknown) point $x\in\mathbb R$ in space at time intervals $(i/N)_{i=0,\dots , N}$ for $N\in\mathbb N$. Then an estimetor
$$\hat H=-\frac{\log 4+\log \left(\frac{1}{N}\sum_{i=0}^{N-1}\left( u\left(\frac{i+1}{N},x\right)-u\left(\frac{i}{N},x\right)\right)^2\right)}{2\log N}$$
of $H$ is strongly consistent and asymptoticaly normal for $\frac{1}{2} \leq H<\frac{3}{4}$. Moreover, it is also strongly consistent for $H>\frac{3}{4}$.
\end{theorem}
\begin{figure}[H]
\centering
\includegraphics[width=0.45\textwidth]{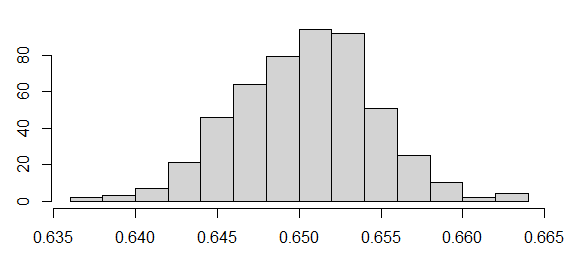}
\includegraphics[width=0.45\textwidth]{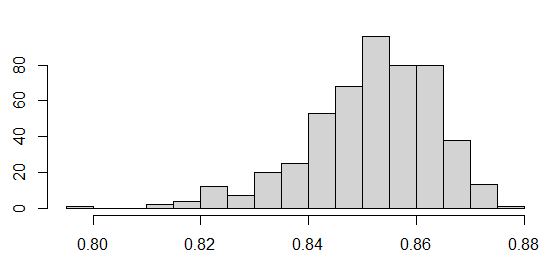}
\caption{Histograms for $500$ independent copies of $\hat H$ for $H=0.65$ and $0.85$ respectively.}
\label{fig1}
\end{figure}
\begin{figure}[H]
\centering
\includegraphics[width=0.45\textwidth]{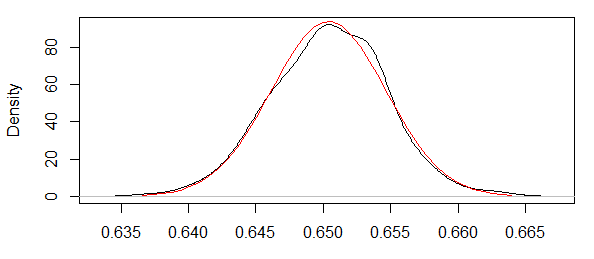}
\includegraphics[width=0.45\textwidth]{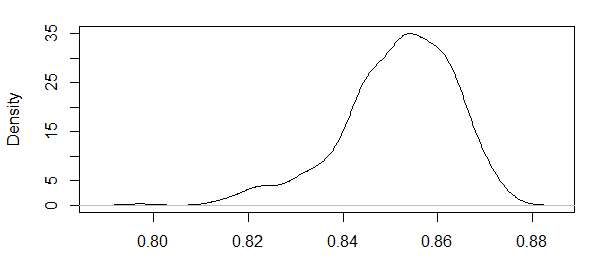}
\caption{Empirical densities of $\hat H$ for $H=0.65$ (with a fitted normal density) and $0.85$ respectively.}
\label{fig2}
\end{figure}
Figures \ref{fig1} and \ref{fig2} obtained with simulations in \verb|R| for $N=1000$ illustrate the performance of the estimator described in the above theorem, in particular they provide a visualisation for different asymptotic behaviour of $\hat H$ for $H$ on either side of the threshold $\frac{3}{4}$.
\begin{remark}
Recall that the limiting variance of the estimator in the paper \cite{KT} was equal to
$$\frac{1}{2}\sum_{k\in\mathbb Z}\varphi_H(k)^2.$$
Thus, for all values of $H$ the estimator based on time increments has a lower asymptotic variance.
\end{remark} 
\section{Estimation of a linear factor}\label{EstimC}
In this section we consider the wave equation with a linear drift parameter and show how quadratic variations can be used for its estimation.We restrict ourselves to the case of a known $H\in \Big[\frac{1}{2},\,\frac{3}{4}\Big)$.

The equation
\begin{equation}
\left\{
\begin{array}{rcl}\label{systeme wave c}
\frac{\partial^2 u}{\partial t^2}(t,x)&=&c^2 \frac{\partial^2 u}{\partial x^2}(t,x)
+\dot W^H(t,x),\;t\geq 0,\;x \in \mathbb{R},\nonumber\\
\noalign{\vskip 2mm}
u(0, x)&=& 0, \quad x \in \mathbb{R},\nonumber\\
\noalign{\vskip 2mm} \frac{\partial u}{\partial t}(0,x) &=& 0,\quad
x \in \mathbb{R}
\end{array} \right.
\end{equation}
for $c\neq 0$  has the solution
$$u^c(x,t)=\frac{1}{2c}\int_0^t\int_{\mathbb R}1_{\{|y-x|\leq c(t-s)\}}dW^H (y,s).$$
This process is well defined because it has the same distribution as the process $u^1=u$ rescaled appropriately:
\begin{equation}\label{eq:sim}
(u^c(x,t))_{t\geq 0, x\in\mathbb R}\stackrel{d}{=}(c^{-1\slash 2}u(\frac{x}{c}, t))_{t\geq 0, x\in\mathbb R}.
\end{equation}
One can see this by using the self-similarity property of the space component of $W^H$:
\begin{eqnarray*}
&& (u^c(x,t))_{t\geq 0, x\in\mathbb R}=\left(\frac{1}{2c}\int_0^t\int_{\mathbb R}1_{\{|y-x|\leq c(t-s)\}}dW^H (y,s)\right)_{t\geq 0, x\in\mathbb R}\\
&&\qquad \left(\frac{1}{2c}\int_0^t\int_{\mathbb R}1_{\{|\frac y c -\frac x c |\leq (t-s)\}}dW^H (y,s)\right)_{t\geq 0, x\in\mathbb R}\\
&&\qquad \stackrel{d}{=} \left(\frac{1}{2c}\int_0^t\int_{\mathbb R}1_{\{|y -\frac x c |\leq (t-s)\}}\sqrt{c} dW^H (y,s)\right)_{t\geq 0, x\in\mathbb R}=(c^{-1\slash 2}u(\frac{x}{c}, t))_{t\geq 0, x\in\mathbb R}.
\end{eqnarray*}
The equality in distribution can be concluded by comparing the covariances, since the processes on both sides are Gaussian. Indeed, the (deterministic) change of variables formula for $\frac{y}{c}\mapsto y$ creates an additional factor $c$ in the covariance integral.

We define
\[V_N^c:=\frac{1}{N}\sum_{i=0}^{N-1}\left(u^c\left(\frac{i+1}{N},\,x\right)-u^c\left(\frac{i}{N},\,x\right) \right)^2-\mathbb E\left[\left(u^c\left(\frac{i+1}{N},\,x\right)-u^c\left(\frac{i}{N},\,x\right) \right)^2\right]\]
as well as
$$S^c_N:=\frac{1}{N}\sum_{i=0}^{N-1}\left( u^c\left(\frac{i+1}{N},x\right)-u^c\left(\frac{i}{N},x\right)\right)^2.$$
Because of the relation \eqref{eq:sim} we have
$$\mathbb E[S^c_N]=c^{-1}\mathbb E [S_N] = c^{-1}\left(\frac{1}{4}N^{-2H}+\frac{1}{4}\frac{1-2H}{1+2H}N^{-2H-1}\right).$$
Therefore, we can define an estimator
$$\hat c:=\frac{1}{4}\frac{1}{S_N^c N^{2H}}.$$
It is strongly consistent, since $V_N^c N^{2H}\to 0$ almost surely by the same argument as above, $V_N^c$ being a rescaled version of $V_N$. Moreover, for $H<\frac{3}{4}$ we have
\[cV_N^c N^{2H+\frac{1}{2}}\stackrel{d}{\to} N(0,\sigma^2),\]
following the same rescaling argument. From this we can obtain the asymptotic normality property for the estimator $\hat c$. We calculate
\begin{eqnarray*}
&&\sqrt N (\hat c-c)=\sqrt{N}\frac{1}{4}\left(\frac{1}{(V_N^c+\mathbb E S_N^c)N^{2H}}-\frac{1}{\mathbb E S_N^c N^{2H}}\left(1+\frac{1-2H}{1+2H}N^{-1}\right)\right)\\
&&\qquad \sqrt{N}\frac{1}{4N^{2H}}\left(\frac{V_N^c}{(V_N^c+\mathbb E[S_N^c])\mathbb E [S_N^c]}-\frac{1-2H}{1+2H} \frac{1}{\mathbb E [S_N^c]N} \right).
\end{eqnarray*}
The second summand tends to zero as $N$ goes to infinity, and thus $\sqrt N (\hat c-c)$ has the same limit as
\begin{eqnarray*}
&& \sqrt{N}\frac{1}{4N^{2H}} V_N^c\frac{1}{(V_N^c+\mathbb E[S_N^c]) c^{-1}\left(\frac{1}{4N^{2H}}+\frac{1-2H}{1+2H}\frac{1}{4N^{2H+1}}\right)}\\
&&\qquad = \sqrt{N} V_N^c c\frac{1}{(V_N^c+\mathbb E[S_N^c]) \left(1+\frac{1-2H}{1+2H}\frac{1}{N}\right)}\\
&&\qquad = \sqrt{N} V_N^c c N^{2H}\frac{1}{(N^{2H}V_N^c+N^{2H}\mathbb E[S_N^c]) \left(1+\frac{1-2H}{1+2H}\frac{1}{N}\right)}.
\end{eqnarray*}
The fraction
$$\frac{1}{(N^{2H}V_N^c+N^{2H}\mathbb E[S_N^c]) \left(1+\frac{1-2H}{1+2H}\frac{1}{N}\right)}$$
converges almost surely to $4c$, and hence, by Slutsky's lemma, we obtain
$$\sqrt N (\hat c-c)\stackrel{d}{\to}N(0,16c^2\sigma^2).$$
In total, we obtain the following theorem.
\begin{theorem}
Consider the wave equation \eqref{systeme wave c} with a known Hurst parameter $H$ in the interval $\Big[\frac{1}{2},\,\frac{3}{4}\Big)$ and an unknown parameter $c\neq 0$. Assume that one realisation of its solution $\Big(u^c\Big(\frac{i}{N},\,x\Big)\Big)_{i=1,\dots ,N}$ for $N\in\mathbb N$ for some possibly unknown $x\in\mathbb R$ is observed. Then
$$\hat c=\frac{1}{4}\frac{1}{ N^{2H-1}\sum_{i=0}^{N-1}\left( u^c\left(\frac{i+1}{N},x\right)-u^c\left(\frac{i}{N},x\right)\right)^2}$$
is a strongly consistent, asymptotically normal estimator of $c$.
\end{theorem}
\begin{remark}
	Following the same scheme we can obtain a consistent, asymptotically normal estimator for the parameter $p:=c\Sigma^{-2}$ in the equation
	\begin{equation}
	\left\{
	\begin{array}{rcl}\label{systeme wave cSigma}
	\frac{\partial^2 u}{\partial t^2}(t,x)&=&c^2 \frac{\partial^2 u}{\partial x^2}(t,x)
+\Sigma\dot W^H(t,x),\;t\geq 0,\;x \in \mathbb{R},\nonumber\\
	\noalign{\vskip 2mm}
	u(0, x)&=& 0, \quad x \in \mathbb{R},\nonumber\\
	\noalign{\vskip 2mm} \frac{\partial u}{\partial t}(0,x) &=& 0,\quad
	x \in \mathbb{R}
	\end{array} \right.
	\end{equation}
	with the solution
	$$u^{c,\Sigma}(x,t)=\frac{\Sigma}{2c}\int_0^t\int_{\mathbb R}1_{\{|y-x|\leq c(t-s)\}}dW^H (y,s),$$
assuming, again, that $H\in \Big[\frac{1}{2},\,\frac{3}{4}\Big)$ is known and we have equidistant temporal observations in the $[0,1]$ interval at one single point in space.

	Moreover, since we have
	$$\mathbb P^{(u^{c,\Sigma}(x,t))_{t\in\mathbb R^+}}=\mathbb P^{(u^{\bar c,\bar \Sigma}(x,t))_{t\in\mathbb R^+}}$$
	for any $x\in\mathbb R$ if and only if $\Sigma c^{-1/2}= \bar\Sigma \bar c^{-1/2}$, there is no hope of estimating these parameters separately. This is similar to the observations made in \cite{HT} for the white-white heat equation. However, with one set each of temporal and spatial observations one can, indeed, recover both quantities: For a fixed time $t$ the quantity
	$$\hat{q}:=\frac{t}{2}\frac{1}{N^{2H}\frac{1}{N}\sum_{i=0}^{N-1}\left(u(t, \frac{i+1}{N})- u(t, \frac{i}{N})\right)^2}$$
	defines a strongly consistent, asymptotically normal estimator for $q:=c^2\Sigma^{-2}$. This can be shown completely analogously by using results from \cite{KT} using the fact that
	$$\mathbb E \left[ \frac{1}{N}\sum_{i=0}^{N-1}\left(u(t, \frac{i+1}{N})- u(t, \frac{i}{N})\right)^2 \right] \sim \frac{t}{2}\frac{\sigma^2}{c^2}\frac{1}{N^{2H}}.$$
\end{remark}
\section{Rectangular quadratic variations: white noise case}\label{S7}
From now on we consider the white noise case, i.e. $H=\frac 1 2$. In this simpler framework we will prove a central limit theorem for rectangular increments with some restrictions for the possible discretisation schemes. As before, we start with an analysis of the covariance structure. For simplicity, we denote in this section $W:=W^{\frac{1}{2}}$.
\begin{lemma}\label{CovRect}
	The solution field $(u(t,\,x))_{t\in\mathbb R^+, \, x\in\mathbb R}$ of the white-white wave equation has the covariance structure
\[\mathbb E[u(t,\,x)u(s,\,y)]= \frac{s^2}{4}\]
for $t\geq s$ and $t-s\geq|x-y|$ and
\[\mathbb E[u(t,\,x)u(s,\,y)]= \frac{(t+s-|x-y|)^2}{16}1_{\{|x-y|<t+s\}}\]
for $t\geq s$ and $t-s\ < |x-y|$.
\end{lemma}

\begin{proof}
Recall that
$$\mathbb E[u(t,\,x)u(s,\,y)]=\int_0^{t\wedge s}\int_\mathbb R G_1(t-u, x-z)G_1(s-u, y-z)dzdu= \frac{1}{4}\int_0^{t\wedge s}\int_\mathbb R 1_{\{|x-z|<t-u\}}1_{\{|y-z|<s-u\}}dzdu.$$
Let us assume without loss of generality that $s\leq t$. The proof is a straightforward calculation where we consider different cases.\\
Case 1: $x\geq y$, $t-s>x-y$. Then the integral becomes
$$\frac{1}{4}\int_0^s \int_{u-s+y}^{s-u+y}dzdu=\frac{s^2}{4}.$$
Case 2: $x\geq y$, $t-s<x-y$.If $x-y<s+t$, the integral equals
$$\frac{1}{4}\int_0^{\frac{s+t+y-x}{2}}\int_{u-t+x}^{s-u+y}dzdu=\frac{(s+y+t-x)^2}{16},$$
otherwise it is zero.\\
Case 3:  $x< y$, $t-s>x-y$. This integral is identical to the one in case 1.\\
Case 4: $x< y$, $t-s<x-y$.Similarly to case 3, if $y-x<s+t$, the integral equals
$$\frac{1}{4}\int_0^{\frac{s+t+x-y}{2}}\int_{u-s+y}^{t-u+x}dzdu=\frac{(s+t+x-y)^2}{16},$$
otherwise it is zero.
\end{proof}

Note that considering $u(x,\cdot)$ for a fixed space point $x$ yields a process equivalent in distribution to a time-rescaled Brownian motion, namely $(B_{t^2\slash 2})_{t\geq 0}$ for a Brownian motion $B$. Since its increments are independent, one can obtain the CLT proved in Section \ref{CLT} by applying the classical CLT and rescaling appropriately. For the estimator of the parameter $c$ defined above we obtain in this special case
\[\sqrt N (\hat c-c)\to N(0,\,\frac{8}{3}c^2).\]
For the quadratic variations of the process $u(\cdot ,t)$ with a fixed time $t$ a CLT has been proved in \cite{KTZ}. More precisely, it has been demonstrated that for
\[V_N=\sum_{j=0}^{N-1}\left[ \frac{(u(t,\,\frac{j+1}{N})-u(t,\,\frac{j}{N}))^2}{\mathbb E[(u(t,\,\frac{j+1}{N})-u(t,\,\frac{j}{N}))^2]}-1\right]\]
one has
\[d\left(\frac{1}{\sqrt{2N}}V_N, N(0,1)\right)\lesssim \frac{1}{\sqrt{N}}.\]
Similarly to Section \ref{EstimC}, one can translate this result to the solution of the equation \eqref{systeme wave c} and subsequently construct an estimator for $c$, utilising the fact that
\begin{eqnarray*}
&&\mathbb E \left[\sum_{j=0}^{N-1}(u^c(t,\,\frac{j+1}{N})-u^c(t,\,\frac{j}{N}))^2\right]\\
&&\qquad=c^{-1}\mathbb E \left[\sum_{j=0}^{N-1}(u(t,\,\frac{j+1}{Nc})-u(t,\,\frac{j}{Nc}))^2\right]\stackrel{N\to\infty}{\sim}\frac{t}{2c^2}-\frac{1}{8N^2 c^3}\stackrel{N\to\infty}{\sim}\frac{t}{2c^2}.
\end{eqnarray*}
It can be checked by direct calculation that the estimator defined as \[\sqrt{\frac{t}{2\sum_{j=0}^{N-1}(u^c(t,\,\frac{j+1}{N})-u^c(t,\,\frac{j}{N}))^2}}\] satisfies the asymptotic normality property with the rate of convergence $\sqrt N$.

Let us now formulate and prove an exemplary CLT for rectangular increments. For simplicity we consider the setting with observations on the grid
$$\mathcal G:=\left\{\left(j\delta_M,\, i\Delta_N\right),\, i=1,\dots , N,\, j=1,\dots ,M\right\},$$
such that for any $(t,x),\,(s,y)\in\mathcal G$ we have $t+s<|x-y|$ if $x\neq y$, that is, for $M=M(N)\stackrel{N\to\infty}{\to}\infty$, we assume $2\delta_M M <\Delta_N$. For instance, this is satisfied for $M=N$, $\Delta_N=\frac{1}{N}$ and $\delta_M=\frac{1}{4N^2}$. In the following we assume another (non-restrictive, but handy) simplification and let $\Delta_N=\frac{1}{N}$ and $\delta_M=M^{-\alpha}$ for some admissible $\alpha$. Studying the behaviour of variations depending on the mesh size in a more general setup is certainly a promising problem, and we consider this result as a motivation for this line of research.

We consider the random variables
\begin{eqnarray*}
\mathcal V_{MN}:=&&\frac{1}{NM}\sum_{i=0}^{N-1}\sum_{j=0}^{M-1}(u_{x_{i+1}t_{j+1}} - u_{x_{i+1}t_{j}} - u_{x_{i}t_{j+1}} + u_{x_{i}t_{j}})^2\\
&&\quad-\mathbb E[(u_{x_{i+1}t_{j+1}} - u_{x_{i+1}t_{j}} - u_{x_{i}t_{j+1}} + u_{x_{i}t_{j}})^2]
\end{eqnarray*}
with $x_i=i\Delta_N$, $t_j=j\delta_j$ and $u_{xt}=u(t,x)$. By Lemma \ref{CovRect} we have
\begin{equation*}
\mathbb E [\Delta_{(ij)}\Delta_{(kl)}]=\begin{cases} \frac{2j+1}{8M^{2\alpha}} \text{ if }i=k,\,j=l, \\ -\frac{2j+1}{8M^{2\alpha}}\text{ if }i\in \{k-1,k+1\}, j=l,\\ 0 \text{ else}\end{cases}
\end{equation*}
with the notation
\[\Delta_{(ij)}:=u_{x_{i+1}t_{j+1}} - u_{x_{i+1}t_{j}} - u_{x_{i}t_{j+1}} + u_{x_{i}t_{j}}.\]
For the asymptotic variance of $\mathcal V_{MN}$ we thus obtain by Gaussianity
\begin{eqnarray*}
&&\mathbb E [\mathcal V_{MN}^2]=\frac{1}{N^2M^2}\sum_{i,k=0}^{N-1}\sum_{j,l=0}^{M-1}2\mathbb E [\Delta_{(ij)}\Delta_{(kl)}]^2\\
&&\qquad \stackrel{N\to\infty}{\sim} \frac{2}{N^2M^2} \sum_{j=0}^{M-1}\sum_{i=1}^{N-1}(\mathbb E [\Delta_{(ij)}^2]^2+ \mathbb E [\Delta_{(ij)}\Delta_{(i+1j)}]^2+ \mathbb E [\Delta_{(ij)}\Delta_{(i-1j)}]^2)\\
&&\qquad =\frac{2}{N^2M^2}\sum_{j=0}^{M-1} N\cdot \frac{(2j+1)^2}{M^{4\alpha}}\frac{3}{64}=\frac{3}{32}\frac{1}{N^2 M^{2+4\alpha}}\sum_{j=0}^{N-1}(2j+1)^2\\
&&\qquad \stackrel{N\to\infty}{\sim}\frac{3}{32}\frac{1}{NM^{2+4\alpha}}\frac{4M^3}{3}=\frac{1}{8}N^{-1}M^{1-4\alpha}.
\end{eqnarray*}
Let us set $\mathcal F_{MN}:=\sqrt{8N}M^{2\alpha-1/2} \mathcal V_{MN}$ to be the normalised version of $\mathcal V_{MN}$. We can now demonstrate the announced CLT.
\begin{theorem}\label{CLTRectangle}
We have under the above assumptions for the solution $u$ of the wave equation driven by a space-time white noise $d(\mathcal F_{MN},N(0,1))\lesssim \frac{1}{\sqrt{NM}}$.
\end{theorem}
\begin{proof}
We note that
\begin{eqnarray*}
&&\Delta_{(ij)}=\int_{\mathbb R^2} G_1(t_{j+1}-r, x_{i+1}-y)1_{r\in [0, t_{j+1}]}dW(r,y)\\
&&\quad -\int_{\mathbb R^2} G_1(t_{j+1}-r, x_{i}-y)1_{r\in [0, t_{j+1}]}dW(r,y) - \int_{\mathbb R^2} G_1(t_{j}-r, x_{i+1}-y)1_{r\in [0, t_{j}]}dW(r,y)\\
&&\quad+ \int_{\mathbb R^2} G_1(t_{j}-r, x_{i}-y)1_{r\in [0, t_{j}]}dW(r,y)\\
&&\qquad =:\int_{\mathbb R^2}g_{ij}(r,y)dW(r,y)
\end{eqnarray*}
with $g_{i,j}$ defined accordingly. Let us denote by $\mathbb H$ the Hilbert space on which the isonormal Gaussian family $W$ is defined. Moreover, let us write $I^W_q$ ($q\in\mathbb N$) for multiple Wiener-It\-o integrals with respect to $W$. Then we can write, similarly to the proof of Theorem \ref{clt-time},
\[\mathcal F_{MN}=N^{-1/2}M^{2\alpha -3/2} \sum_{i=0}^{N-1}\sum_{j=0}^{N-1}I^W_1(g_{ij})^2-\mathbb E[I^W_1(g_{ij})^2]=N^2\sum_{i,j=0}^{N-1}I^W_2(g_{ij}\otimes g_{ij})\]
as well as
\[D\mathcal F_{MN}=N^{-1/2}M^{2\alpha -3/2}\sum_{i=0}^{N-1}\sum_{j=0}^{M-1} I^W_1(g_{ij})g_{ij},\]
and consequently,
\[\|D\mathcal F_{MN}\|^2_{\mathbb H}-\mathbb E [\|D\mathcal F_N\|^2_{\mathbb H}]=N^{-1}M^{4\alpha-3} \sum_{i,k=0}^{N-1}\sum_{j,l=0}^{M-1}\langle g_{ij}, g_{kl} \rangle_{\mathbb H} I^W_2(g_{ij}\otimes g_{kl})  \]
as well as
\[\operatorname{Var}(\|D\mathcal F_{MN}\|^2_{\mathbb H})=N^{-2}M^{8\alpha-6}\sum_{i,k,m,o=0}^{N-1}\sum_{j,l,n,p=0}^{M-1} \langle g_{ij}, g_{kl} \rangle_{\mathbb H} \langle g_{mn}, g_{op} \rangle_{\mathbb H} \langle g_{ij}, g_{mn} \rangle_{\mathbb H} \langle g_{kl}, g_{op} \rangle_{\mathbb H}. \]
We know that
\begin{eqnarray*}
&&|\langle g_{ij}, g_{kl} \rangle_{\mathbb H}|=\left| \mathbb E \left[\int_{\mathbb R^2}g_{ij}(r,y)dW(r,y)\int_{\mathbb R^2}g_{kl}(r,y)dW(r,y) \right]\right|\\
&&\qquad =|\mathbb E [\Delta_{(ij)}\Delta_{(kl)}]|\leq \frac{2j+1}{8M^{2\alpha}}1_{\{j=l\}}1_{\{|i-k|\leq 1\}},
\end{eqnarray*}
and thus,
\begin{eqnarray*}
&&\operatorname{Var}(\|D\mathcal F_{MN}\|^2_{\mathbb H}) \lesssim N^{-2}M^{8\alpha-6} \sum_{i=0}^{N-1}\sum_{j=0}^{M-1}\left(\frac{2j+1}{M^{2\alpha}}\right)^4\\
&&\qquad= N^{-2}M^{8\alpha-6} NM^{-8\alpha}\sum_{j=0}^{M-1}(2j+1)^4\lesssim N^{-1}M^{-1}.
\end{eqnarray*}
Consequently,  $d(\mathcal F_{MN},N(0,1))\lesssim \frac{1}{\sqrt{MN}}$.
\end{proof}

In particular, the theorem states that the rescaled quadratic variations on $\mathcal G$ converge to $N(0,1)$ with the same speed as the rescaled quadratic variations on either the time or the space part of the grid. Moreover, this CLT permits us the study of an estimator for the parameter $c$ in the setting \eqref{systeme wave c} defined via rectangular increments. Note first that for
\[\mathcal S_{MN}^c :=\frac{1}{MN}\sum_{i=0}^{N-1}\sum_{j=0}^{M-1}(\Delta^c_{ij})^2\]
with
\[\Delta^c_{(ij)}:=u^c_{x_{i+1}t_{j+1}} - u^c_{x_{i+1}t_{j}} - u^c_{x_{i}t_{j+1}} + u^c_{x_{i}t_{j}}\]
we have
\[ \mathbb E[8M^{2\alpha-1} \mathcal S_{MN}^c ]=8M^{2\alpha-2}N^{-1}\sum_{i=0}^{N-1}\sum_{j=0}^{M-1} \mathbb E[(\Delta^c_{ij})^2]=8c^{-1}M^{2\alpha-2} \sum_{j=0}^{M-1}\frac{2j+1}{8M^{2\alpha}}= c^{-1}(1+O(1/M))\]
if $\frac{\Delta_N}{c}>2\delta_M M$. Here the restriction on $c$ can be weakened by enlarging the spatial observation interval. For such $c$ we can define the estimator
\[\tilde{c}:=\frac{1}{8 M^{2\alpha-1} \mathcal S_{MN}^c }.\]
We obtain pointwise
\begin{eqnarray*}
&&\tilde c-c \sim \frac{1}{8 M^{2\alpha-1}\mathcal S_{MN}^c }- \frac{1}{8 M^{2\alpha-1}\mathbb E [\mathcal S_{MN}^c ]}\\
&&=\frac{1}{8 M^{2\alpha-1}}\left(\frac{\mathbb E [\mathcal S_{MN}^c]- \mathcal S_{MN}^c}{\mathcal S_{MN}^c \mathbb E [\mathcal S_{MN}^c] }\right)=-\frac{\mathcal V_{MN}^c}{8 M^{2\alpha-1} \mathbb E [\mathcal S^c_{MN}](\mathcal V_{MN}^c + \mathbb E [\mathcal S_{MN}^c])},
\end{eqnarray*}
where
\[\mathcal V^c_{MN}=\frac{1}{MN}\sum_{i=0}^{N-1}\sum_{j=0}^{N-1}(\Delta^c_{ij})^2-\mathbb E [(\Delta^c_{ij})^2].\]
Recalling that by hypercontractivity \eqref{hyper}
\[ M^{2\alpha-1} \mathcal V^c_{MN}\to 0 \text{ a.s.}\]
and by the CLT (Theorem \ref{CLTRectangle})
\[c\sqrt{8N}M^{2\alpha-1/2}\mathcal V^c_{MN}\stackrel{d}{\to}N(0,1),\]
we conclude
\[\sqrt{MN}(\tilde c-c)\sim -\frac{\sqrt{MN}\mathcal V_{MN}^c}{8 M^{2\alpha-1} \mathbb E [\mathcal S^c_{MN}](\mathcal V_{MN}^c + \mathbb E [\mathcal S_{MN}^c])} \stackrel{d}{\to}N(0,\,8c^2).\]
These statistical results are summarised below.
\begin{theorem}
Consider the wave equation \eqref{systeme wave c} driven by white spatiotemporal noise, i.e. with $H=\frac{1}{2}$, and an unknown drift parameter $c\neq 0$. Suppose we observe its solution $u^c$ on a space time grid
$$\mathcal G=\left\{ (t_j, x_i) = \left(jM^{-\alpha},\, iN^{-1}\right),\, i=1,\dots , N,\, j=1,\dots ,M\right\},$$
such that $M(N)\to\infty$ when $N\to\infty$ and $2M^{1-\alpha}<N^{-1}$. Then for $N\to\infty$ the estimator
\[\tilde{c}=\frac{1}{8 M^{2\alpha-2}N^{-1}\sum_{i=0}^{N-1}\sum_{j=0}^{M-1}(\Delta^c_{ij})^2 }\]
with
\[\Delta^c_{(ij)}=u^c(t_{j+1}, x_{i+1}) - u^c(t_{j}, x_{i+1}) - u^c(t_{j+1}, x_{i}) + u^c(t_{j}, x_{i})\]
is a strongly consistent and asymptotically normal estimator of the drift parameter $c$.
\end{theorem}
The limiting variance of the resulting estimator is larger than that for the estimator obtained from the observations in time. However, this result demonstrates that one can indeed use a mix of space and time observations to assess the drift parameter.

\section{Appendix}

The basic tools from the analysis on Wiener space are presented in this section. We will focus on some elementary facts about multiple stochastic integrals. We refer to \cite{N} for a complete review on the topic. 

Consider ${\mathcal{H}}$ a real separable infinite-dimensional Hilbert space
with its associated inner product ${\langle
.,.\rangle}_{\mathcal{H}}$, and $(B (\varphi),
\varphi\in{\mathcal{H}})$ an isonormal Gaussian process on a
probability space $(\Omega, {\mathfrak{F}}, \mathbb{P})$, which is a
centered Gaussian family of random variables such that
$\mathbb{E}\left( B(\varphi) B(\psi) \right) = {\langle\varphi,
\psi\rangle}_{{\mathcal{H}}}$ for every
$\varphi,\psi\in{\mathcal{H}}$. Denote by $I_{q}$ the $q$th multiple
stochastic integral with respect to $B$, which is an
isometry between the Hilbert space ${\mathcal{H}}^{\odot q}$
(symmetric tensor product) equipped with the scaled norm
$\frac{1}{\sqrt{q!}}\Vert\cdot\Vert_{{\mathcal{H}}^{\otimes q}}$ and
the Wiener chaos of order $q$, which is defined as the closed linear
span of the random variables $H_{q}(B(\varphi))$ where
$\varphi\in{\mathcal{H}},\;\Vert\varphi\Vert_{{\mathcal{H}}}=1$ and
$H_{q}$ is the Hermite polynomial of degree $q\geq 1$ defined
by:\begin{equation}\label{Hermite-poly}
H_{q}(x)=(-1)^{q} \exp \left( \frac{x^{2}}{2} \right) \frac{{\mathrm{d}}^{q}%
}{{\mathrm{d}x}^{q}}\left( \exp \left(
-\frac{x^{2}}{2}\right)\right),\;x\in \mathbb{R}.
\end{equation}The isometry of multiple integrals can be written as follows: for $p,\;q\geq
1$,\;$f\in{{\mathcal{H}}^{\otimes p}}$ and
$g\in{{\mathcal{H}}^{\otimes q}}$
\begin{equation} \mathbb{E}\Big(I_{p}(f) I_{q}(g) \Big)= \left\{
\begin{array}{rcl}\label{iso}
q! \langle \tilde{f},\tilde{g}
\rangle _{{\mathcal{H}}^{\otimes q}}&&\mbox{if}\;p=q,\\
\noalign{\vskip 2mm} 0 \quad\quad&&\mbox{otherwise}.
\end{array}\right.
\end{equation}It also holds that:
\begin{equation*}
I_{q}(f) = I_{q}\big( \tilde{f}\big),
\end{equation*}
where $\tilde{f} $ denotes the canonical symmetrization of $f$ and is defined by $$\tilde{f}%
(x_{1}, \ldots , x_{q}) =\frac{1}{q!}\sum_{\sigma\in\mathcal{S}_q}
f(x_{\sigma (1) },\ldots, x_{\sigma (q)}),$$ where the sum runs
over all permutations $\sigma$ of $\{1,\ldots,q\}$.

We have the following product formula: if  $f\in{{\mathcal{H}}^{\odot p}}$ and
$g\in{{\mathcal{H}}^{\odot q}}$, then 

\begin{eqnarray}\label{prod}
I_{p}(f) I_{q}(g)&=& \sum_{r=0}^{p \wedge q} r! \binom{p}{r}\binom{q}{r}I_{p+q-2r}\left(f\tilde{\otimes}_{r}g\right)
\end{eqnarray}

A further important  property of  finite sums of multiple integrals is the hypercontractivity. Namely, if $F= \sum_{k=0} ^{n} I_{k}(f_{k}) $ with $f_{k}\in \mathcal{H} ^{\otimes k}$ then
\begin{equation}
\label{hyper}
\mathbb{E}\vert F \vert ^{p} \leq C_{p} \left( \mathbb{E}F ^{2} \right) ^{\frac{p}{2}}.
\end{equation}
for every $p\geq 2$.

We denote by $D$ the Malliavin derivative operator that acts on
cylindrical random variables of the form $F=g(B(\varphi
_{1}),\ldots,B(\varphi_{n}))$, where $n\geq 1$,
$g:\mathbb{R}^n\rightarrow\mathbb{R}$ is a smooth function with
compact support and $\varphi_{i} \in {{\mathcal{H}}}$. This
derivative is an element of $L^2(\Omega,{\mathcal{H}})$ and it is
defined as
\begin{equation*}
DF=\sum_{i=1}^{n}\frac{\partial g}{\partial x_{i}}(B(\varphi _{1}),
\ldots , B(\varphi_{n}))\varphi_{i}.
\end{equation*}

\noindent
\textbf{Acknowledgements}
The author would like to thank Dr. Meryem Slaoui for helpful discussions and the two anonymous referees for the useful corrections and suggestions. The support of the German Research Foundation (DFG) via SFB 823 is gratefully acknowledged.

%% Bibliography

%%  For BibTeX users:
 \bibliographystyle{unsrt}
 \bibliography{biblio}

%\bibliography{biblio}{}
%\bibliographystyle{alpha} 

\end{document}